\documentclass[12pt]{article}
\input amssymb.sty
\input amssym.def
\input amssym
\input epsf

\newtheorem{theorem}{Theorem}[section]
\newtheorem{lemma}[theorem]{Lemma}
\newtheorem{proposition}[theorem]{Proposition}
\newtheorem{corollary}[theorem]{Corollary}
\newtheorem{conjecture}[theorem]{Conjecture}
\newtheorem{definition}[theorem]{Definition}

\begin{document}
\newcommand{\be}{\begin{equation}}
\newcommand{\ee}{\end{equation}}
\newcommand{\bt}{\begin{theorem}}
\newcommand{\et}{\end{theorem}}
\newcommand{\bd}{\begin{definition}}
\newcommand{\ed}{\end{definition}}
\newcommand{\bp}{\begin{proposition}}
\newcommand{\ep}{\end{proposition}}
\newcommand{\bl}{\begin{lemma}}
\newcommand{\el}{\end{lemma}}
\newcommand{\bc}{\begin{corollary}}
\newcommand{\ec}{\end{corollary}}
\newcommand{\bcon}{\begin{conjecture}}
\newcommand{\econ}{\end{conjecture}}
\newcommand{\la}{\label}
\newcommand{\Z}{{\Bbb Z}}
\newcommand{\R}{{\Bbb R}}
\newcommand{\Q}{{\Bbb Q}}
\newcommand{\C}{{\Bbb C}}
\newcommand{\hra}{\hookrightarrow}
\newcommand{\lra}{\longrightarrow}
\newcommand{\lms}{\longmapsto}
\newcommand{\AAA}{{\Bbb A}}
\newcommand{\PP}{{\Bbb P}}

\begin{titlepage}
\title{  A simple construction of Grassmannian polylogarithms}
\author{A.B. Goncharov }

\date{\it  To Andrey Suslin for his 60th birthday}
 
\end{titlepage}
\maketitle

\tableofcontents

\begin{abstract}
The classical $n$-logarithm is a multivalued analytic 
function defined inductively:
$$
{\rm Li}_n(z):= \int_0^z {\rm Li}_{n-1}(t)d\log t, \quad {\rm Li}_1(z)= - \log (1-z).
$$
In this paper we give a simple explicit construction of the Grassmannian 
$n$-logarithm, which is a multivalued 
 analytic function on   the quotient of the Grassmannian of 
$n$-dimensional subspaces in  $\C^{2n}$ in generic position 
to the coordinate hyperplanes by the natural action of the  torus 
$(\C^*)^{2n}$. The classical $n$-logarithm 
appears at a certain one dimensional boundary stratum. 
\vskip 2mm

We study {\it Tate iterated integrals}, which are 
homotopy invariant integrals of $1$-forms $d\log f_i$ where $f_i$ are rational functions. 
We give a simple explicit formula for the Tate iterated integral 
which describes the Grassmannian 
$n$-logarithm. 

Another example is the Tate iterated integrals 
for the 
multiple polylogarithms 
on the moduli spaces ${\cal M}_{0,n}$, 
calculated in Section 4.4 of \cite{G2} using the combinatorics 
of  plane trivalent trees 
decorated by the arguments of the  multiple polylogarithms. 

\vskip 2mm 
Variations of mixed Hodge-Tate structures on $X$ are described by a Hopf algebra 
${\cal H}_\bullet(X)$. 
We upgrade Tate iterated integrals on a (rational) complex variety $X$ to elements of  ${\cal H}_\bullet(X)$. 
The coproducts of these elements are very interesting invariants of the iterated integrals. 
In general their calculation is a nontrivial problem.  
We show however, that working modulo the ideal of ${\cal H}_\bullet(X)$ generated by constant variations, 
 there is a simple way to calculate the coproduct. 

\vskip 2mm
It is a pleasure to 
dedicate this paper to Andrey Suslin, 
whose works \cite{Su} and \cite{Su2} played an essential role 
in the development of the story. 
\end{abstract}
\section{Introduction and main definitions}

\subsection{The Grassmannian polylogarithms and their properties}

\paragraph{Configurations and Grassmannians.} 
A {\it configuration} of $m$ points of a $G$-set $X$ is an
orbit of the group $G$ on $X^m$. Recall the classical  
dictionary relating configurations of points in projective/vector spaces to 
Grassmannians.

1. If $X= V_n$ is an $n$-dimensional complex vector space and $G= GL_n(\C)$ 
we have 
configurations  of vectors in $V_n$. Configurations of vectors in isomorphic 
vector spaces are canonically identified. 
Such a  configuration is {\it generic} if any $k\leq n$ vectors are linearly independent. 

Denote by $G_{n}$ the moduli space of generic 
configurations of $2n$ vectors in an $n$-dimensional vector space, with respect to the group $GL_n$. 
Its complex points are identified with the points of the open part of the Grassmannian $G^n_n(\C)$
of $n$-dimensional subspaces in the coordinate space $\C^{2n}$ parametrising the subspaces 
which are in generic position to the coordinate hyperplanes. 
Namely, such a 
 subspace $H \subset \C^{2n}$ 
provides a configuration of $2n$ vectors in the dual space 
$H^*$ given by the restriction of the coordinate functions. 
\vskip 2mm

2. If $X= \C\PP^{n-1}$, $n>1$, and $G= PGL_n(\C)$ we have 
configurations  of points in $\C\PP^{n-1}$. 
Such a  configuration is {\it generic} if any $k \leq n$ of
the points generate a $(k-1)$-plane in $\C\PP^{n-1} $. 

Denote by $P G_{n}$ the moduli space of generic 
configurations of $2n$ points in $\PP^{n-1}$. 
Its complex points are identified with the orbits of the torus $({\C}^*)^{2n}$ acting on the Grassmannian 
$G^n_n(\C)$. Namely, an $n$-dimensional
 subspace $H \subset \C^{2n}$ 
provides a configuration of $2n$ hyperplanes in the projectivisation of 
$H$ given by intersection with the coordinate hyperplanes. 
By the projective duality this is the same as a generic configuration of $2n$ points in 
$\C\PP^{n-1}$.

\paragraph{Construction of the Grassmannian polylogarithms.} 
The Grassmannian $n$-logarithm is a multivalued 
analytic function $L_n^G$ on $P G_{n}(\C)$, which 
we define as the integral of an explicit closed $1$-form $\Omega$ on $PG_n(\C)$. 
The $1$-form $\Omega$ is defined by 
using the Aomoto $(n-1)$-logarithms \cite{A}, whose definition we recall now.

\paragraph{\it The Aomoto $n$-logarithm.} 
A simplex $L$ in $\C{\Bbb P}^{n}$ is a collection of $n+1$ hyperplanes $(L_0, ..., L_{n})$. 
In particular, a collection of $n+1$ points in generic position 
determines a simplex with the vertices at these points. 
A pair of simplices $(L; M)$ in $\C{\Bbb P}^{n}$ is {\it admissible} if $L$ and $M$
have no common faces of the same dimension.
 There is a canonical $n$-form $\omega_L$ in ${\Bbb
C}{\Bbb P}^n$ with logarithmic poles at the hyperplanes $L_i$. Namely, if $z_i = 0$ are
homogeneous equations of $L_i$ then 
$$
\omega_L = d\log(z_1/z_0) \wedge
...  \wedge d\log(z_n/z_0).
$$
 Recall that for a nondegenerate simplex $M$, the rank of the relative homology group 
${\rm rk}H_n({\Bbb C}{\Bbb P}^n,  M )$  is one. 
Let $\Delta_M$ be a topological $n$-cycle representing
a generator of $H_n({\Bbb C}{\Bbb P}^n,  M)$. 
The Aomoto $n$-logarithm is a
multivalued analytic function on configurations of admissible pairs of simplices
$(L;M)$ in ${\Bbb C}{\Bbb P}^n$ given by 
$$
{\cal A}_n(L;M) := \int_{\Delta_M} \omega_L. 
$$

\vskip 2mm
{\bf Examples}. 1. Let $(l_1,l_2)$ and $(m_1,m_2)$ be two pairs of distinct 
points in $\C{\Bbb P}^1$. 
Then 
$$
{\cal A}_1(l_1,l_2; m_1,m_2) := \int_{m_1}^{m_2} d \log \frac{z-l_2}{z-l_1} = \log r(l_1,l_2, m_1,m_2).
$$
where $r(x_1,x_2, x_3,x_4)$ is the cross-ratio of four points on the projective line: 
$$
r(x_1,x_2, x_3,x_4):= \frac{(x_3-x_1)( x_4-x_2)}{(x_3-x_2)
(x_4-x_1)}. 
$$
Here $\int_{m_1}^{m_2}$ denotes the integral along a path connecting $m_1$ and $m_2$, 
which does not contain the other two points.

2. The classical $n$-logarithm ${\rm Li}_n(z)$ is given by an $n$-dimensional integral
$$
{\rm Li}_n(z) = \int_{ 0  \leq 1 - t_1 \leq t_2 \leq ... \leq t_n \leq z}
\frac{d t_1}{t_1}\wedge ... \wedge\frac{d t_n}{t_n}. 
$$

\vskip 2mm

Below we always use the following convention about the integration cycles $\Delta_M$. 
Given a generic configuration of points $(x_1, ..., x_m)$ in $\C{\Bbb P}^{n-1}$, 
a {\it compatible system of chains} is the following data.  
For every two points $(x, y)$ of the configuration 
we choose a generic oriented path $\varphi(x,  y)$ connecting them, 
for every three points $(x, y, z)$ we choose 
a generic oriented topological triangle $\varphi(x,  y, z)$ 
which bounds $\varphi(x,  y) + \varphi(y, z) + \varphi(z, x)$, and so on,  
so that for every subconfiguration $(x_{i_1}, ..., x_{i_k})$, $k\leq n$ we choose a 
generic oriented topological simplex 
$\varphi(x_{i_1}, ..., x_{i_k})$, 
and these choices are compatible with the boundaries. 
In the definition of the Aomoto polylogarithms we always 
choose a $\varphi$-simplex as the chain $\Delta_M$. 

\vskip 2mm 

Let $V_n$ be an $n$-dimensional complex 
vector space. 
Choose a volume form $\omega_n \in {\rm det} ~V_n^*$. 
Given vectors $l_1,...,l_n$ in $V_n$, set 
$$
\Delta(l_1,...,l_n):= \langle l_1\wedge ... \wedge l_n, \omega_n \rangle. 
$$
Notice that a simplex in a projective space ${\Bbb P}(V)$ 
can be defined as either a collection of hyperplanes, 
or vertices. Below we employ the second point of view, and use vectors $l_i\in V$ 
to determine the vertices 
as the lines spanned by the vectors. 
 
Consider the following multivalued analytic 
$1$-form on the Grassmannian $G_n(\C)$: 
\begin{equation} \label{gras1}
\Omega(l_1,...,l_{2n}):= {\rm Alt}_{2n} \Bigl(
{\cal A}_{n-1}(l_1,...,l_{n}; l_{n+1}, ... , l_{2n}) ~d\log 
\Delta(l_{n+1},...,l_{2n})\Bigr).
\end {equation}
Here ${\rm Alt}_{2n}$ denotes the alternation of a function in $2n$ variables, that is the 
alternated sum of $(2n)!$ terms. 
It does not depend on the choice of the form $\omega_n$, since the latter 
does not vary, and appears under the differential.

\begin {theorem} \label{MMTH} 
For any $l_1, ..., l_{2n}$ in generic position in an $n$-dimensional complex vector space, 
the $1$-form  $\Omega(l_1,...,l_{2n})$ is closed. 
It depends only on the configuration of points in $\C{\Bbb P}^{n-1}$ 
obtained by projection of the vectors $l_i$. 
\end {theorem}


\begin {definition} \label{MMTH1} The Grassmannian $n$-logarithm $L_n^G(l_1,...,l_{2n})$ is 
the  skew-symmetrization under 
the permutations of the vectors $l_1, ..., l_{2n}$ of 
the primitive of the $1$-form (\ref{gras1}). 
\end {definition}

A primitive of the  $1$-form (\ref{gras1}) is a multivalued analytic function 
defined up to a scalar. 
The scalar vanishes under the skew-symmetrization. So 
the Grassmannian $n$-logarithm is a well defined multivalued analytic function. 
Thanks to the last claim of Theorem \ref{MMTH} we can consider it as a function $L_n^G(x_1,...,x_{2n})$ 
on configurations of $2n$ points in $\C {\Bbb P}^{n-1}$.

\paragraph{Properties of the Grassmannian $n$-logarithm.}

Given a configuration of $m+1$ vectors $(l_0, ..., l_m)$ in $V_n$, denote by 
$(l_0|l_1, ..., l_m)$ a configuration of vectors obtained by projection of the 
vectors $l_1, ..., l_m$ to the quotient of $V_n$ along the subspace 
generated by $l_0$. We employ a projective version of this construction. 
Given a configuration of $m+1$ points $(y_0, y_1,...,y_m)$ in 
$\C{\Bbb P}^{n-1}$, denote by 
$(y_0 \vert y_1,...,y_m)$ the  configuration of $m$ points in $\C{\Bbb P}^{n-2}$ 
obtained by projection of the points $y_i$ with the center at the point 
$y_0$.

\begin {theorem} \label{MMTH2} 

The function $L_n^G(x_1,...,x_{2n})$ enjoys the following properties.

\begin{enumerate} 

\item {\rm The $(2n+1)$-term equation.} For a generic 
configuration of $2n+1$ points $(x_1,...,x_{2n+1})$ in 
$\C\PP^{n-1}$ one has 
$$
\sum_{i=1}^{2n+1}(-1)^i L_n^G(x_1,...,\widehat x_i, ... , x_{2n+1}) = \mbox{ \rm a constant}.
$$
\item {\rm Dual $(2n+1)$-term equation.} For a  generic configuration 
of points $(y_1,...,y_{2n+1})$ in $\C\PP^n$ 
$$
\sum_{j=1}^{2n+1}(-1)^j L_n^G(y_j| y_1,...,\widehat y_j, ... , y_{2n+1}) = \mbox{ \rm a constant}.
$$
\end{enumerate}
\end {theorem}

Here we assumed that compatible systems of cycles for the configurations 
of points $(x_1,...,x_{2n+1})$ and $(y_1,...,y_{2n+1})$ were chosen. 

{\bf Example}. For $n=2$ we get the Rogers version of the dilogarithm:
$$
L_2^G(x_1,x_2,x_3,x_4) = L_2(r(x_1,x_2,x_3,x_4)),\quad 
\mbox{where}\quad  L_2(z):= {\rm Li}_2(z) + \frac{1}{2}\log(1-z)\log(z).
$$

\begin{figure}[ht]
\centerline{\epsfbox{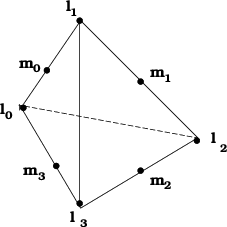}}
\caption{Special configuration of $8$ points in ${\Bbb P}^3$.}
\label{fig1coc}
\end{figure}

A configuration $(x_1,...,x_n,y_1,...,y_n)$ of points 
in ${\Bbb P}^{n-1}$  is called a {\it special configuration} if 
$(x_1, ..., x_n)$ form a 
generic configuration, and for every $i$ the point 
$y_i$ lies on the line $x_i ~x_{i+1}$, where the indices are modulo $n$. 
 See an example on Fig \ref{fig1coc}. Special configurations are parametrised  
by one  parameter, denoted by  $r(x_1,...,x_n,y_1,...,y_n)$, see \cite{G4}, Section 4.4.
For $n=2$ it is the cross-ratio. 
One can show ({\it loc. cit}) that the restriction of the function $L_n^G$ to a special 
configuration is 
expressed via the classical $n$-logarithm function.

The Grassmannian $n$-logarithm is a period of a variation of framed mixed $\Q$-Hodge-Tate 
 structures of geometric origin 
on $PG_{n}(\C)$. 
We call it {\it the Grassmannian variation of mixed Tate motives}. 
Below we introduce and  calculate the Tate 
iterated integral related to the Grassmannian polylogarithm function.

\subsection{The history and ramifications of the problem.} 
There are three incarnations of the dilogarithm function: 

\vskip 2mm
i) The real valued Rogers dilogarithm ${\rm L}_2(x)$ defined 
on $\R{\Bbb P}^1 - \{0, 1, \infty\}$ by the condition:
\be \la{DE}
d{\rm L}_2(x) = \frac{1}{2}\Bigl(-\log|1-x|d\log |x|+ \log|x|d\log |1-x|\Bigr), \qquad 
{\rm L}_2(-1) = {\rm L}_2(1/2) = {\rm L}_2(2) =0. 
\ee
Notice that $\R{\Bbb P}^1 - \{0, 1, \infty\}$ is 
the moduli space of generic configuration of $4$ points in $\R{\Bbb P}^1$, for the group $PGL_2(\R)$.  
The function ${\rm L}_2(r(l_1, ..., l_4))$ is the unique solution 
of the differential equation (\ref{DE}) which is skew symmetric 
under the permutations of the vectors $l_i$. 
Its restriction to the interval $(0,1)$ is given by 
$$
{\rm L}_2(x) = {\rm Li}_2(x) + \frac{1}{2}\log (1-x)\log x - \frac{\pi^2}{12}, \qquad x \in (0,1).
$$
 It satisfies the $5$-term relation
$$
\sum_{k=0}^4(-1)^k 
{\rm L}_2(r(l_0, ..., \widehat {l_k}, ..., l_4)) = 
-\varepsilon \frac{\pi^2}{6}, \qquad \varepsilon = 
\frac{1}{2}\prod_{0 \leq i<j \leq 4}{\rm sgn} ~\Delta(l_i, l_j).
$$ 
 \vskip 2mm
ii) The multivalued complex analytic 
dilogarithm function ${\rm Li}_2(z)$, whose properties are 
best described by the corresponding variation of framed mixed $\Q$-Hodge structures.

\vskip 2mm
iii) The single valued Bloch-Wigner function, defined on $\C{\Bbb P}^1$ by 
$$
{\cal L}_2(z) := {\rm Im}\Bigl(Li_2(z) + \log(1-z) \log|z|\Bigr).
$$
It satisfies the $5$-term relation
$$
\sum_{k=0}^4(-1)^k 
{\cal L}_2(r(l_0, ..., \widehat {l_k}, ..., l_4)) = 0.
$$ 
The Bloch-Wigner function is nothing else but the real period of the variation 
which appears in ii). 
The motivic nature of the dilogarithm is described by the 
Bloch-Suslin complex and its relations to  algebraic K-theory,  see \cite{Su2}. 

\vskip 3mm

In accordance with this, there are three directions  
for a generalization of the dilogarithm function: 

\vskip 2mm
i) Gelfand and MacPherson [GM] defined a  
real valued Grassmannian $2n$-logarithm function 
on $PG_{2n}(\R)$ by constructing its differential. 
Notice that our construction of the Grassmannian $n$-logarithm 
also starts from a closed $1$-form $\Omega$ on $G_n(\C)$. 
The relationship between these two functions is not clear. It should reflect the relationship 
between the Chern and Pontryagin classes.

\vskip 2mm
ii) The construction of Hanamura and MacPherson \cite{HM1}, \cite{HM2} provides a
Grassmannian $n$-logarithm function. 
The  construction is geometric but rather complicated. 
I do not know how to relate it to the function $L_n^G$. 
An explicit motivic construction of a Grassmannian $n$-logarithm function 
was given for $n=3$ in \cite{G} and for $n=4$ in \cite{G3}. 

\vskip 2mm
iii) In \cite{G1}, see also \cite{G4}, we defined 
a single-valued Grassmannian $n$-logarithm function  ${\rm L}^G_{n, \R}$. 
The precise relationship between this function 
and the multivalued analytic 
function ${\rm L}_n^G$  is not known. 

\vskip 3mm
\paragraph{The bi-Grassmannian $n$-logarithm cocycles.} We denote by $G_p^q$ the 
Grassmannian of $q$-dimensional subspaces in a coordinate vector space of 
dimension $p+q$, transversal to the coordinate hyperplanes. 
The weight $n$ bi-Grassmannian ${\Bbb G}(n)_\bullet^\bullet$ is 
given by a collection of Grassmannians $G_p^q$, 
$p \geq n$, arranged in a form of a truncated bisimplicial variety:
\be \la{BG}
\begin{array}{ccccccccc}
&&& &\ldots &&\ldots&&\ldots\\
&&&&&&\downarrow &&\downarrow \\
&&G_{n+1}^{n}&\lra& \ldots &\lra&G_{n+1}^{2}&\lra&G_{n+1}^{1}\\
&&\downarrow &&\ldots &&\downarrow &&\downarrow \\
G_n^{n+1}&\lra&G_n^{n}&\lra &\ldots &\lra &G_n^{2} &\lra & G_n^{1}\\
\end{array}
\ee
Here a horizontal arrow stands for a collection of maps 
given by the intersection of the subspaces with the coordinate hyperplanes, and the vertical one for 
projection along the coordinate axes, see \cite{G5}. 
The bottom line is the semi-simplicial weight $n$ Grassmannian ${\Bbb G}^\bullet_n$ introduced in \cite{BMS}. 

The weight $n$ bi-Grassmannian ${\Bbb G}(n)_\bullet^\bullet$ 
and the related polylogarithms play a key role in 
the explicit combinatorial construction of Chern classes suggested in \cite{G5}.

\vskip 2mm
Points of the bi-Grassmannian (\ref{BG}) with values in 
a field $F$ form a 
truncated bisimplicial set. Applying to it the ``free abelian group'' functor $S \to \Z[S]$ 
we get a bi-Grassmannian complex. 
Its bottom line is the Grassmannian complex, whose homology was studied by Suslin in the 
fundamental paper \cite{Su}. 
\vskip 2mm

Each of the three versions of the Grassmannian $n$-logarithm functions 
should appear as a component of the corresponding 
 {\it bi-Grassmannian $n$-logarithm cocycle}, which is a 
cocycle in the complex calculating the cohomology 
of the bi-Grassmannian with coefficients in a certain complex of sheaves. 
These complexes are: 

\vskip 2mm
i) A real analog of the weight $2n$ Deligne complex on ${\Bbb G}(n)_\bullet^\bullet(\R)$.

ii) The multivalued analytic weight $n$ Deligne complex on ${\Bbb G}(n)_\bullet^\bullet(\C)$ considered in \cite{BMS}.

iii) The real weight $n$ Deligne complex on ${\Bbb G}(n)_\bullet^\bullet(\C)$ -- see, for example, \cite{G1}.
\vskip 2mm

Here is what is known about the corresponding cocycles.

\vskip 2mm
i) The real bi-Grassmannian $2n$-logarithm cocycle is the crucial building block 
in the Gabrielov, Gelfand and Losik \cite{GGL} approach to 
a combinatorial formula for the $n$-th Pontryagin class. 
However such a cocycle 
is available only when $2n=2$ \cite{GGL}, and, mostly, when $2n=4$ 
\cite{Yu}, \cite{G3}. 

\vskip 2mm
ii) The existence of a multivalued analytic 
Grassmannian $n$-logarithm cocycle was conjectured by 
 Beilinson, MacPherson and Schechtman  \cite{BMS}. 
An explicit geometric construction was found 
in \cite{HM1}, 
\cite{HM2}. A weaker existence theorem was proved in \cite{H}. 
There is an explicit motivic construction of the bi-Grassmannian $n$-logarithm cocycle
for $n=3$ \cite{G} and $n=4$ \cite{G3}.

\vskip 2mm
iii) A single-valued bi-Grassmannian $n$-logarithm cocycle was defined in \cite{G1}, see also \cite{G4}. 
It has a rather peculiar property: its components assigned to the 
Grassmannians $G_m^\bullet(\C)$, $m>n$ (i.e. above the bottom row in 
(\ref{BG}))   are identically  zero. 
This is not expected to hold for the motivic/multivalued  analytic  bi-Grassmannian 
$n$-logarithm cocycles for $n>3$.

\subsection{Symbols and Tate iterated integrals}
In Section 3 we introduce {\it Tate iterated integrals} on a complex algebraic variety $X$. They  
are certain (conjecturally all) homotopy invariant iterated 
integrals of $1$-forms $d\log f_i$, where 
$f_i$ are rational functions on $X$. 

Denote by ${\cal O}_X^*$  
the multiplicative group of regular invertible functions on $X$.  
The length $n$ Tate iterated integrals are 
determined by their {\it symbols}
\be \la{EII}
{\rm I} \in \bigotimes^n \overline {{\cal O}_X^*}, \qquad \overline {{\cal O}_X^*}:= {\cal O}_X^*/\C^*,
\ee
satisfying certain integrability condition of algebraic nature. For $n=2$ 
the integrability just means 
that the image of the element ${\rm I}$ in $K_2(X)$ modulo the symbols $\{\C^*, {\cal O}_X^*\}$ in $K_2$  
is zero, see Definition 3.1. 

Beilinson's construction (cf. \cite{DG}) implies that any 
Tate iterated integral is the  period of an $n$-framed variation ${\cal I}({\rm I})$ 
of mixed  motives on $X \times X$, 
understood as a variation of mixed $\Q$-Hodge structures of geometric origin. 

We show that when $X$ is rational, there is   an $n$-framed geometric  variation  
of mixed $\Q$-Hodge-Tate structures whose period is the Tate iterated integral.  
Conjecturally the same is true for any $X$, justifying the  name. 

Conversely, any variation ${\cal V}$ of $n$-framed 
mixed $\Q$-Hodge-Tate structures on complex manifold $M$ determines a symbol 
\be \la{EIIa}
{S}_n({\cal V}) \in 
\bigotimes^n{\cal O}^*_{M, {\rm an}}. 
\ee
Here ${\cal O}^*_{M, {\rm an}}$ is the multiplicative group of invertible analytic functions on $M$. 

The targets of the symbols (\ref{EII}) and (\ref{EIIa}) are different: 
in (\ref{EII}) we kill the constants, while in (\ref{EIIa}) we do not, and the symbol 
(\ref{EIIa}) is an analytic object. Moreover,  although an analytic  
symbol ${\Bbb I} \in \bigotimes^n{\cal O}^*_{M, {\rm an}}$ produces 
an iterated integral  on $M$, in general it is not a period of a variation 
of mixed Hodge structures. 

\vskip 2mm
So there are  two constructions: 

\begin{itemize}

\item A symbol ${\rm I}$ on $X$ provides a variation on $X \times X$ with the fiber 
${\cal I}_{x,y}({\rm I})$ at $(x,y) \in X\times X$; 

\item We assign 
to a variation on $X$ a symbol (\ref{EIIa}) on $X$. 

\end{itemize}

They are related as follows. 
Given a point $a$, there is a geometric Hodge-Tate variation ${\cal I}_{a,y}({\rm I})$ on $X$. 
Considered modulo the ideal generated by constant variations on $X$, it does not depend on $a$ (Lemma 3.8). 
The symbol ${\cal S}({\cal I}({\rm I}))$ lies in  $\bigotimes^n{\cal O}_X^*$, its projection to 
$\bigotimes^n\overline {{\cal O}_X^*}$ does not depend on $a$ and equals to the original symbol ${\rm I}$. 

\vskip 3mm
Framed variations of $\Q$-Hodge-Tate structures on $X$ give rise to a commutative 
graded Hopf algebra ${\cal H}_\bullet(X)$. Any $n$-framed $\Q$-Hodge-Tate variation on $X$ 
provides an element of ${\cal H}_n(X)$. 
In general the coproduct of 
the element ${\cal I}({\rm I}) \in {\cal H}_n(X)$ corresponding to an integrable symbol ${\rm I}$ is 
rather complicated.  
We show, however, that, considered modulo the ideal generated by the constant variantions, 
the coproduct is determined by the deconcatenation map on the symbols -- see Theorem \ref{VITJ}.  

The symbol of the Tate iterated integral correponding to the period of a
 geometric variation of $\Q$-Hodge-Tate variation on $X$ 
can be calculated inductively if we know the differential equation of the period function. 

\vskip 3mm
{\bf Conclusion}. Working modulo the ideal of the Hopf algebra ${\cal H}_\bullet(X)$ generated by the constant variations, 
 we arrive at a simple and effective way to calculate the coproducts 
of the elements of  ${\cal H}_\bullet(X)$ corresponding to periods of geometric Hodge-Tate variations. 

\paragraph{The structure of the paper.} In Section 2 we 
recall the  scissors congruence groups ${A}_n(F)$, whose properties reflect  
the ones of the Aomoto $n$-logarithm. 
The functional equations of the Grassmannian $n$-logarithm stated in Theorem \ref{MMTH2} 
follow immediately from basic properties of the Aomoto $(n-1)$-logarithm.  
However Theorem \ref{MMTH}, and therefore the existence 
of the Grassmannian $n$-logarithm function $L_n^G$,  is less obvious.  
It is proved in Section 2. 

 In Section 3 we discuss symbols and {\it Tate iterated integrals}.  

In Section 4 we 
define explicitly 
a Tate iterated integral on the Grassmannian $G_{n}(\C)$ by exhibiting its 
symbol ${\rm I}_n$. 
We prove that  
${\rm I}_n$ coincides with the symbol
 of the iterated integral provided by the integration of the form $\Omega$.

\paragraph{Acknowledgments.} 
I was 
supported by the  NSF grants DMS-0653721 and DMS-1059129. 
This paper was written at the IHES (Bures sur Yvette)  
during the Summer of 2009. 
I am grateful to IHES  for the support. 
I am very much indebted to the referee for many useful comments, which improved the exposition.

\section{Properties of the Grassmannian polylogarithms}

\subsection{Motivic avatar of the form $\Omega$} 

\paragraph{The scissors 
congruence groups $A_n(F)$.} They were defined in 
\cite{BMS}, \cite{BGSV}. We use slightly modified groups, 
adding one more relation -- the dual additivity relation.  

Let $F$ be a field. 
The abelian group $A_n(F)$ is generated by the elements $$
\langle l_0,...,l_n;m_0,...,m_n\rangle _{A_n}
$$ corresponding to generic  
configurations of $2(n+1)$ points
$(l_0,...,l_n;m_0,...,m_n)$ in ${\Bbb P}^n(F)$. We use the notation
  $\langle L;M\rangle _{A_n}$ where $L= (l_0,...,l_n)$ and $M = (m_0,...,m_n)$. 
The relations, which reflect properties of the Aomoto polylogarithms, are 
the following:

\begin{enumerate}

\item 
{\it Nondegeneracy}.  $\langle L;M\rangle _{A_n} = 0$ if $(l_0,...,l_n)$ or $(m_0,...,m_n)$ belong to a
hyperplane.

\item 
{\it Skew symmetry}. $\langle \sigma L;M\rangle _{A_n} = \langle L;\sigma M\rangle _{A_n} = (-1)^{|\sigma|} \langle L;M\rangle _{A_n}$
for any $\sigma \in S_{n+1}$.

\item 
{\it Additivity}.  For any configuration $(l_0,...,l_{n+1})$ 
$$
\sum_{i=0}^{n+1}(-1)^i \langle l_0,...,\widehat l_i,...,l_{n+1};m_0,...,m_n\rangle _{A_n} =0,
$$ 
and a similar
condition for $(m_0,...,m_{n+1})$. 

{\it Dual additivity}. For any configuration $(l_0,...,l_{n+1})$ 
$$
\sum_{i=0}^{n+1}(-1)^i \langle l_i\vert l_0,...,\widehat l_i,...,l_{n+1};m_0,...,m_n\rangle _{A_n} =0,
$$
and a similar
condition  for $(m_0,...,m_{n+1})$.

\item 
{\it Projective invariance}. $\langle gL;gM\rangle_{A_n}  = \langle L;M\rangle_{A_n}$ for any $g \in PGL_{n+1}(F)$. 
\end{enumerate}

The cross-ratio provides a
canonical isomorphism
$$
a_1: A_1(F) \longrightarrow F^{\ast}, \quad a_1: \langle l_0,l_1;m_0,m_1\rangle _{A_1} \longmapsto 
r(l_0,l_1,m_0,m_1).
$$

\bl
The Aomoto polylogarithm function satisfies all the above properties 1)-4).
\el

{\bf Proof}. Follows straight from the definitions. Notice that it is essential to use 
the compatible system of topological simplices $\varphi$ as 
representatives of the relative cycles $\Delta_M$. 

\paragraph{The coalgebra $A_\bullet(F)$.} Set $A_0(F) = \Bbb Z$. 
There is a graded coassociative 
coalgebra structure on  $A_\bullet(F):= \oplus_{n \geq 0} A_n(F)$ 
with a coproduct $\nu$, see 
\cite{BMS}, \cite{BGSV}.\footnote{The coproduct $\nu$ is defined by the same formula as 
in {\it loc. cit.}. Recall that the formula works only for generic pairs of simplices. 
The combinatorial formula for the coproduct used in {\it loc. cit.} 
in the Hodge or $l$-adic realizations coincides with (and was motivated by) the 
general formula for the coproduct of framed objects in mixed categories, \cite{G2}, Appendix. 
The derivation of the former from the latter is a good exercise. A detailed solution of a similar 
problem for a different kind of 
 scissor congruence groups is given in Theorem 4.8 in \cite{G6}. 
See Section 4 there for further details.} 
We need only one component of the coproduct:
$$
\nu_{n-1, 1} : A_n(F) \lra A_{n-1}(F)\otimes_\Z F^*.
$$
We employ a formula for 
$\nu_{n-1, 1}$ derived in Proposition 2.3 of \cite{G3}, 
which is much more convenient than the original one for computations and manifestly skew-symmetric. 
Using the notation ${\rm Alt}_{3,3}$ for the skew-symmetrization 
of $(l_0, l_1, l_2)$ as well as $(m_0, m_1, m_2)$, we have\footnote{The coefficient $-1/4$ in (\ref{BNU}) is 
compatible with the specialisation of formula (\ref{n>2coproduct}) for $\nu_{n-1,1}$ plus a similar 
formula for $\nu_{1,n-1}$ for $n=2$. Indeed, if $n=2$, there are 
$3\times 3$ terms in each of the two formulas, total $18$, 
while in  (\ref{BNU}) the total number of terms, before taking $1/4$, is $2(3!)^2=72$.}
\be \la{BNU}
\nu_{1,1} \langle l_0, l_1, l_2; m_0, m_1, m_2\rangle _{A_2} = 
\ee
$$
-\frac{1}{4}{\rm Alt}_{3,3}\Bigl(\Delta(m_0, l_1, l_2) \otimes \langle m_0|l_1, l_2;
m_1, m_2\rangle _{A_1} + 
\langle l_0|l_1, l_2;
m_1, m_2\rangle _{A_1} \otimes \Delta(l_0,m_1, m_2) \Bigr).
$$
For $n>2$:
\be \la{n>2coproduct}
\nu_{n-1,1}\Bigl(\langle l_0,...,l_n;m_0,...,m_n\rangle _{A_n}\Bigr) = 
\ee
$$
-\sum_{i,j =0}^n (-1)^{i+j}\langle l_i\vert l_0,...,\hat 
l_i,...,l_n;m_0,...,\hat m_j,...,m_n\rangle _{A_{n-1}}   \otimes
\Delta(l_i,m_0,...,\hat
m_j,...,m_n  ).
$$
It is straightforward to prove that $\nu_{n-1,1}$ 
is well defined, i.e. kills the relations.

\paragraph{The map $\nu_{n-1, 1}$ and the differential of the Aomoto polylogarithm.}
Let ${\Bbb A}_n$ be the field of rational 
functions on the space of pairs of simplices in $\C{\Bbb P}^n$.  
There is a natural map
$$
{\cal A}_{n} \otimes d\log: {A}_n({\Bbb A}_{n}) \otimes_\Z {\Bbb A}^*_{n} \lra 
\Omega^1_{\rm mv}, \qquad \langle L, M\rangle \otimes F \lms 
{\cal A}_{n} (L, M) ~d\log (F).
$$
where $\Omega^1_{\rm mv}$ is the space of multivalued analytic $1$-forms 
on the space of pairs of simplices in $\C{\Bbb P}^n$.  
\bl \la{6.29.09.3} One has 
\be \la{dfsdf}
d {\cal A}_n(l_0, ..., l_n; m_0, ..., m_n) = {\cal A}_{n-1} \otimes d\log 
~ \circ ~ \nu_{n-1, 1} \langle l_0, ..., l_n; m_0, ..., m_n\rangle_{{A}_{n-1}}.
\ee
\el

{\bf Proof}. This is a very special case of the general formula for 
the differential of the period of a variation of Hodge-Tate structures, 
see Lemma \ref{3}. 

One can  easily prove it directly as follows. We can assume that 
the vectors $l_0, ..., l_n$ form a standard basis. Let us consider a small 
deformation $m_i(t)$ of the vectors $m_i$, where $0 \leq t \leq \varepsilon$. 
By Stokes formula, 
to calculate the differential of the function ${\cal A}_n(l_0, ..., l_n; m_0(t), ..., m_n(t))$ 
we have to calculate the linear in $\varepsilon$ term of $\sum(-1)^j\int_{M_j(\varepsilon)}\omega_L$, where 
$M_j(\varepsilon)$ is the $n$-dimensional body  obtained 
by moving the $j$-th face $(m_0(t), ..., \widehat m_j(t), ..., m_n(t))$. 
One can easily see that this matches the $j$-th term in (\ref{dfsdf}). 
The lemma is proved. 

\paragraph{Motivic avatar of the form $\Omega$.} Recall the notation  $\Q(X)$ 
for the field of rational functions on a variety $X$ over $\Q$. 
Consider the 
following element of 
\begin{equation} \label{gras111}
\Lambda_{n-1,1}(l_1,...,l_{2n}) \in {A}_{n-1}(\Q(G_n))
\otimes_\Z \Q(G_n)^*.
\end {equation}
\begin{equation} \label{gras11}
\Lambda_{n-1,1}(l_1,...,l_{2n}):= {\rm Alt}_{2n} \Bigl(\langle l_1,...,l_{n}; 
l_{n+1}, ... , l_{2n}\rangle_{A_{n-1}} \otimes \Delta(l_{n+1},...,l_{2n})\Bigr).
\end {equation}
 
\begin{lemma} \label{keyp1}
For any $2n+1$ vectors $(l_1, ..., l_{2n+1})$ in generic position in $V_{n}$ one has 
$$
\sum_{i=1}^{2n+1}(-1)^i\Lambda_{n-1,1} (l_1, ... , \widehat l_i, ... , l_{2n+1}) = 0.
$$
For any $2n+1$ vectors $(m_1, ..., m_{2n+1})$ in generic position in $V_{n+1}$ one has 
$$
\sum_{j=1}^{2n+1}(-1)^j\Lambda_{n-1,1} (m_j| m_1, ... , \widehat m_j, ... , m_{2n+1}) = 0.
$$
\end{lemma}

{\bf Proof}. The first formula follows immediately from the statement 
that 
$$
\sum_{i=1}^{n+1}(-1)^i\langle l_1, ... , \widehat l_i,  ... , l_{n+1}; l_{n+2}, ... , l_{2n+1}\rangle _{A_{n-1}} 
\otimes \Delta(l_{n+2}, ... , l_{2n+1}) = 0
$$
which follows from the additivity. 
The second reduces to the dual additivity. The lemma is proved.

\subsection{Proof of Theorems \ref{MMTH} and \ref{MMTH2}}
Below we always work modulo $2$-torsion. 

\vskip 2mm
We start from the following observations. 
Let ${\cal A}$ be a coassociative  coalgebra with the coproduct $\nu$, and 
${\cal A}_+$  the kernel of the counit. Let 
$$
\widetilde \nu:= \nu - ({\rm Id} \otimes 1 + 1 \otimes {\rm Id}): {\cal A}_+ \lra {\cal A}_+^{\otimes 2}
$$
 be the restricted coproduct. Then there is a map 
$\nu_{[k]}: {\cal A}_+ \lra \otimes^k {\cal A}_+$ 
given by a composition 
$$
{\cal A}_+ \stackrel{\widetilde \nu}{\lra} {\cal A}_+\otimes {\cal A}_+ \stackrel{\widetilde \nu \otimes{\rm Id}}{\lra} {\cal A}_+\otimes {\cal A}_+ \otimes {\cal A}_+ \stackrel{\widetilde \nu \otimes {\rm Id}}{\lra} ... \stackrel{\widetilde \nu \otimes {\rm Id}}{\lra} {\cal A}_+^{\otimes k}.
$$
The coassociativity of ${\cal A}$ implies that
one can replace anywhere here ${\widetilde \nu \otimes {\rm Id}}$ by ${{\rm Id}\otimes \widetilde \nu }$. 

In particular, if
 ${\cal A}:= \oplus {\cal A}_n$ is graded by positive integers, there is  a map  
(\cite{G2}):
$$
\nu_{[n]}: {\cal A}_n \lra \otimes^n {\cal A}_1.
$$

Given an abelian group $A$, there is a commutative graded Hopf algebra given by the 
tensor algebra ${\rm T}(A)$ of   
$A$ with the shuffle product $\circ$ and the coproduct 
$\delta$ given by the  deconcatenation map 
$$
\delta: a_1 \otimes \ldots \otimes a_n \lms \sum_{k=0}^n a_1 \otimes \ldots \otimes a_k \bigotimes a_{k+1} \otimes \ldots \otimes a_n. 
$$ 

\bl \la{11.1.10.5}
Let ${\cal A}:= \oplus_{n\geq 0} {\cal A}_n$ be a connected commutative graded Hopf algebra. Then the map 
$$
\nu: {\cal A} \lra {\rm T}({\cal A}_1)
$$
 given by the direct sum of the maps $\mu_{[n]}$, is a morphism of graded commutative Hopf algebras. 
\el

{\bf Proof}. The claim that $\nu$ commutes with the coproducts follows from 
the very definition. The claim that $\nu$ commutes with the products is easy to check. 
The Lemma is proved. 

\paragraph{Proof of Theorem \ref{MMTH}.}
The case $n=2$ is trivial. For example, the form is closed since in this case 
we deal with functions of one variable. So 
we assume below $n\geq 3$. 
One has 
$$
(  \nu_{n-2,1} \otimes {\rm Id} ) \circ \Lambda_{n-1,1} (l_1,...,l_n;m_1,..,m_n) = 
$$
\be \la{6.29.09.1}
(  \nu_{n-2,1} \otimes {\rm Id} ) 
{\rm Alt}_{2n}
\Bigl(\langle l_1, ... ,l_n; m_{1}, ... ,m_{n}\rangle _{A_{n-1}}  \otimes \Delta (m_{1}, ... ,m_{n})\Bigr) = 
\ee
$$
-n^2\cdot {\rm Alt}_{2n} 
\Bigl(\langle l_1|l_2, ... ,l_n; m_{2}, ... ,m_{n}\rangle _{A_{n-2}}  \otimes \Delta (l_1, m_{2}, ... ,m_{n})\otimes \Delta (m_{1}, ... ,m_{n})\Bigr). 
$$
So thanks to Lemma \ref{6.29.09.3} we need to prove that 
\be \la{222}
{\rm Alt}_{2n}
\Bigl({\cal A}_{n-2}(l_1| l_2, l_3, ... ,l_n; m_2, m_{3}, ... ,m_{n}) 
d\log \Delta (l_1, m_{2}, ... ,m_{n})\wedge d\log \Delta (m_{1}, ... ,m_{n})\Bigr) =0.  
\ee
We will deduce this from the following Lemma
\bl \la{MMLE}
$$
{\rm Alt}_{2n}
\Bigl(d{\cal A}_{n-2}(l_1| l_2, l_3, ... ,l_n; m_2, m_{3}, ... ,m_{n}) 
\otimes d\log \Delta (l_1, m_{2}, ... ,m_{n})\wedge d\log 
\Delta (m_{1}, ... ,m_{n})\Bigr) =0.  
$$
\el

Lemma \ref{MMLE} implies the first claim of Theorem \ref{MMTH} by the following argument: 
Integrating each of the $1$-forms 
$d{\cal A}_{n-2}(l_1| l_2, l_3, ... ,l_n; m_2, m_{3}, ... ,m_{n})$ we recover 
(\ref{222}) plus a sum
$$
\sum C_{\alpha_1, \alpha_2} d\log\Delta_{\alpha_1}\wedge d\log\Delta_{\alpha_2},
$$
where $\alpha_1 = \{l_1, m_2, ..., m_{n}\}$, 
$\alpha_2 = \{m_1, ..., m_n\}$, and $C_{\alpha_1, \alpha_2}$ 
are the integration constants. It is zero 
since we alternate an expression symmetric in $(m_{n-1}, m_n)$.

\vskip 2mm
{\bf Proof of Lemma \ref{MMLE}}. Using (\ref{6.29.09.1}), one has 
$$
(  \nu_{n-3,1} \otimes{\rm Id} \otimes {\rm Id} ) \circ (  \nu_{n-2,1}  \otimes {\rm Id} ) \circ 
\Lambda_{n-1,1} (l_1,...,l_n;m_1,..,m_n) = 
$$
\be \la{PPP}
n^2 (n-1)^2 \cdot {\rm Alt}_{2n}
\Bigl(\langle l_1, l_2|l_3, ... ,l_n; m_{3}, ... ,m_{n}\rangle _{A_{n-3}}  \otimes 
\ee
$$
\Delta (l_1, l_2, m_{3}, ... ,m_{n})\otimes \Delta (l_1, m_{2}, ... ,m_{n})\otimes \Delta (m_{1}, ... ,m_{n})\Bigr).  
$$
 
It is sufficient to prove the following
\bl \la{MMLEa} The element (\ref{PPP}) has zero projection to 
${A}_{n-3}(\Q(G_n)) \otimes 
{\Q(G_n)}^* \otimes K_2(\Q(G_n))$.  
\el

{\bf Proof}. Set $\delta \{x\}:= (1-x) \wedge x$. 
Let us show that, dividing by $n^2 (n-1)^2$, (\ref{PPP}) is equal to 
\be \la{PPPas}
{\rm Alt}_{2n}
\Bigl(\langle l_1, l_2|l_3, ... ,l_n; m_{3}, ... ,m_{n}\rangle _{A_{n-3}}  \otimes 
\Delta (l_1, l_2, m_{3}, ... ,m_{n})\otimes 
\delta \{r(m_3, ... ,m_{n} | l_1, l_2, m_1, m_{2})\}\Bigr).  
\ee
We use the following formula (\cite{G}, Lemma 2.6), valid only modulo $2$-torsion
\footnote{recall that we work modulo $2$-torsion throughout the paper.}:
\be \la{6.29.09.2}
\delta \{r (v_1, v_2, v_3, v_4)\} = \frac{1}{2}{\rm Alt}_4 \Bigl(\Delta(v_1, v_2) \wedge\Delta(v_1, v_3) \Bigr).  
\ee
We say that a single term in formula (\ref{6.29.09.2}),
 say $\Delta(v_1, v_2) \wedge\Delta(v_1, v_3)$, 
 is obtained by choosing $v_1$ and forgetting $v_4$. 

So the 
product of the last two factors in the expression under the alternation sign in 
(\ref{PPP}) is obtained by 
choosing $m_{2}$ and forgetting $l_{2}$ in 
\be \la{df}
\delta \{r(m_3,  ... ,m_{n} | l_1, l_2, m_1, m_{2})\}.
\ee 

1. Due to skew-symmetry, the term obtained by 
choosing $m_{i}$ and forgetting $l_{j}$, where $i=1,2$ and $j=1,2$, also appears. 
We use a similar argument in 2-4 below. 

2. The term obtained by 
choosing $m_{2}$ and forgetting $m_{1}$ vanishes. This follows by applying the 
additivity relation in the first argument to the configuration 
$$
(l_1, l_2 | m_1, l_3, ..., l_n ; m_3, ..., m_n). 
$$  
Indeed, none of the vectors $m_1, l_3, ..., l_n$ enters the last three factors 
(the second row below) 
of the expression 
$$
{\rm Alt}_{2n}\langle l_1, l_2| l_3, ..., l_n; m_3, ..., m_n\rangle_{A_{n-2}} \otimes 
$$
$$
\Delta(l_1,  l_2, m_3, ..., m_n) \otimes \Delta(l_1, m_2, m_3, ..., m_n) \wedge 
\Delta(l_2, m_2, m_3, ..., m_n). 
$$ 

3. The term obtained by 
choosing $l_{1}$ and forgetting $l_{2}$ vanishes. This follows by applying the 
dual additivity relation in the second argument to the configuration 
$$
(l_1 | l_3, ..., l_n ; {l_2}, m_3, ..., m_n).
$$ 
Indeed,  the 
dual additivity relation provides us the first of the following two equalities:
$$
{\rm Alt}_{2n}\langle l_1, l_2| l_3, ..., l_n; m_3, ..., m_n\rangle_{A_{n-2}} \otimes 
$$
$$
\Delta(l_1, l_2, m_3, ..., m_n) \otimes \Delta(l_1, m_1, m_3, ..., m_n) \wedge 
\Delta(l_1, m_2, m_3, ..., m_n) = 
$$ 
$$
-\sum_{k=3}^n(-1)^k{\rm Alt}_{2n}\langle l_1, m_k| l_3, ..., l_n; l_2, m_3, ..., \widehat m_k, ... , m_n\rangle_{A_{n-2}} \otimes 
$$
$$
\Delta(l_1, l_2, m_3, ..., m_n) \otimes \Delta(l_1, m_1, m_3, ..., m_n) \wedge 
\Delta(l_1, m_2, m_3, ..., m_n)  =0. 
$$
To prove the second equality, notice that the pair $(l_1, m_k)$, where $k\geq 3$,
 enters every four factors of the last expression 
symmetrically, and thus the sum vanishes. 

4. The term obtained by 
choosing $l_{1}$ and forgetting $m_{1}$ vanishes. 
This follows by applying the additivity relation 
for the configuration $$
(l_1, l_2 | {m_1}, l_3, ..., l_n ; m_3, ..., m_n).
$$ 
Indeed,  none of the vectors $m_1, l_3, ..., l_n$ enters  the last three factors 
(the second row below) 
of the expression 
$$
{\rm Alt}_{2n}\langle l_1, l_2| l_3, ..., l_n; m_3, ..., m_n\rangle_{A_{n-2}} \otimes 
$$
$$
\Delta( l_1, l_2, m_3, ..., m_n) \otimes \Delta(l_1, l_2, m_3, ..., m_n) \wedge 
\Delta(l_1, m_2, m_3, ..., m_n). 
$$ 
Lemma \ref{MMLEa}, and hence Lemma \ref{MMLE} and the first claim of 
Theorem \ref{MMTH} are proved. 

\vskip 2mm
The form $\Omega$ does not change if we multiply the vector 
$l_{2n}$ by a constant $a \in \C^*$:  
$$
\Omega(l_1, ...  , al_{2n}) - \Omega(l_1, ...,  l_{2n}) = 
{\rm Alt}_{2n-1}\Bigl({\cal A}_{n-1}(l_{1}, ... , l_{n}; l_{n+1}, ..., l_{2n})\Bigr)\otimes d\log a =0.
$$
Indeed, it is easy to prove using Lemma \ref{6.29.09.3} that 
$
{\rm Alt}_{2n}\Bigl(d{\cal A}_{n-1}(l_{1}, ... , l_{n}; m_1, ..., m_n)\Bigr) = 0.
$ 
This implies the claim, just as above. Theorem \ref{MMTH} is proved. 

\begin {conjecture} \label{MMTH3} $\Lambda_{n-1,1}(l_1,...,l_{2n})$ does not change if one of the vectors $l_i$ is multiplied by $\lambda \in F^*$. So it depends only on the configurations of $2n$ points in ${\Bbb P}^{n-1}$ defined by the vectors $l_i$. 
\end {conjecture}

\paragraph{Proof of Theorem \ref{MMTH2}.} Applying the map 
${\cal A}_{n-1} \otimes d\log$ to the element (\ref{gras111}) we get 
the form $\Omega$. Therefore the proof follows from 
Lemma \ref{keyp1}.

\section{Symbols, Tate iterated integrals, and variations of  
mixed Tate motives}

\subsection{Symbols and Tate iterated integrals}

\paragraph{Iterated integrals of smooth $1$-forms.} Let $M$ be a real manifold. Let $\omega_1$, ..., $\omega_n$ be smooth 
$1$-forms on $M$. Then given a path $\gamma: [0,1] \to M$ 
there is an iterated integral 
\be \la{rt}
\int_\gamma \omega_1 \circ  ... \circ \omega_n:= 
\int_{0 \leq t_1 \leq \ldots \leq t_n \leq 1}\gamma^*\omega_1(t_1) \wedge \ldots \wedge 
\gamma^*\omega_n(t_n).
\ee
Let $({\cal A}^*(M), d)$ be the commutative DG $\R$-algebra  of smooth forms on $M$. 
By linearity an element
$$
{\rm I} \in \bigotimes^n({\cal A}^1(M)[1]) := \underbrace{
{\cal A}^1(M)[1]\otimes \ldots \otimes {\cal A}^1(M)[1]}_{\mbox{$n$ factors}}
$$
gives rise to an iterated integral $\int_{\gamma}({\rm I})$. Here $[1]$ stands for the sdhift of grading by one.

\paragraph{Homotopy invariant iterated integrals.} Denote by ${\rm T}(A)$ the tensor algebra of 
a graded vector space $A$. The bar complex of the 
commutative DG algebra ${\cal A}^*(M)$ is defined as 
 ${\rm T}({\cal A}^*(M)[1])$ equipped with a 
 differential  
$$
D : {\rm T}({\cal A}^*(M)[1]) \lra {\rm T}({\cal A}^*(M)[1]). 
$$
The differential 
is the sum of the de Rham differential $d$ and the maps given by 
the products of the consecutive factors in the tensor product. 
A theorem of K.T. Chen \cite{Ch} tells us that an iterated integral $\int_{\gamma}({\rm I})$
is homotopy invariant, i.e. invariant
 under deformations of the path $\gamma$ 
preserving its endpoints,  if and only if $D({\rm I}) =0$. 

In particular,  a collection of 
closed $1$-forms $\omega^{(s)}_i$ such that 
for every $1 \leq k \leq n-1$ one has 
\be \la{123Chen}
\sum_s \omega^{(s)}_1 \otimes  ... \otimes \omega^{(s)}_{k-1} 
\otimes (\omega^{(s)}_{k}\wedge \omega^{(s)}_{k+1}) \otimes \omega^{(s)}_{k+2} \otimes 
...  \otimes \omega^{(s)}_{n} = 0
\ee
gives rise to a homotopy invariant iterated integral 
$
\sum_s \int_{\gamma}\omega^{(s)}_1 \otimes  ... \otimes \omega^{(s)}_n.
$

\paragraph{Symbols and Tate iterated integrals.} 
Now let $X$ be a complex algebraic variety. 
Our goal is to study iterated integrals 
of $1$-forms $d\log f_i$ where $f_i\in {\cal O}_X^*$ are invertible regular functions on $X$. 
There is an inclusion
$$
d\log: \overline {{\cal O}_X^*} \hra \Omega^1_{\rm log}(X), \qquad \overline {{\cal O}_X^*}:= {\cal O}_X^*/\C^*.
$$
Given a path $\gamma: [0,1] \to X(\C)$ in $X(\C)$ and 
$$
{\Bbb I} = f_1(x) \otimes \ldots \otimes f_n(x) \in \bigotimes^n{\cal O}_X^*
$$
there is 
an iterated integral 
$$
\int_{\gamma}d \log ({\Bbb I}) = \int_\gamma d\log f_1 \circ d\log f_2 \circ \ldots \circ d\log f_n.
$$ 
which evidently depends only on the image ${\rm I}$ of the element ${\Bbb I}$ in $\bigotimes^n\overline {{\cal O}_X^*}$.

The forms $d\log f$ are closed.  So condition 
(\ref{123Chen}) implies the homotopy invariance of the corresponding 
iterated integral. The  map $d\log$ annihilates the Steinberg element $(1-f) \otimes f$. 
Conjecturally the ideal generated by the Steinberg elements and constants is
 the kernel of the map $d\log$. So this is an algebraic condition on the functions $f_i$ which 
implies condition 
(\ref{123Chen}), and which is hypothetically  equivalent to it.

This leads to the following definition. 
Let $F$ be a field. Recall that by Matsumoto's theorem, the group $K_2(F)$ is the quotient of $F^* \otimes F^*$ by the subgroup generated by  
the Steinberg relations $(1-x) \otimes x$, where $x \in F^*-\{1\}$. So there is a projection 
$$
\pi: F^* \otimes F^* \lra K_2(F), \qquad a\otimes b \lms \{a,b\}.
 $$
For $1 \leq k \leq n-1$ there is a map obtained by 
applying $\pi$ to the $k$-th factor $\otimes^2F^*$  in $\otimes^nF^*$: 
$$
\bigotimes^n F^* \lra \bigotimes^{k-1} F^* \otimes K_2(F) \otimes \bigotimes^{n-k-1} F^*, \qquad \pi_k = {\rm Id} \otimes \pi \otimes {\rm Id}.
$$
Since we work with $\overline {\C(X)^*}$ rather then with $\C(X)^*$, we take these maps modulo the 
ideal  generated by $\C^*$ in the tensor algebra of $\C(X)^*$. So we arrive at the projections 
$$
\pi_{k,n}: \bigotimes^n \overline {\C(X)^*} \lra \bigotimes^{k-1} \overline {\C(X)^*} 
\otimes \frac{K_2(\C(X))}{\{\C^*, \C(X)^*\}} \otimes \bigotimes^{n-k-1} \overline {\C(X)^*}.
$$
\bd
An element ${\rm I} \in \bigotimes^n\overline {\C(X)^*}$ is integrable if $\pi_{k,n}({\rm I}) =0$ 
for every  $1 \leq k \leq n-1$.\footnote{In the case $n=1$ any element is integrable.}  
An integrable symbol on $X$ is an element 
\be \la{IEL}
{\rm I} \in \bigotimes^n\overline {{\cal O}_X^*}
\ee
 whose image in $\bigotimes^n\overline {\C(X)^*}$ 
is integrable. 
\ed

\bd A {\rm Tate iterated integral} is an iterated integral 
given by an integrable element (\ref{IEL}).  
The element ${\rm I}$ is called the {\rm symbol} of the Tate iterated integral.  
\ed
Chen's theorem immediately implies that Tate iterated integrals are homotopy invariant.

\subsection{The Hopf algebra of integrable symbols.}

Consider the direct sum  
\be \la{11.1.10.1}
{\bf I}_\bullet(X):= \bigoplus_{n=0}^\infty{\bf I}_n(X), \qquad  
{\bf I}_n(X):= {\rm Int}\Bigl(\bigotimes^n\overline {{\cal O}_X^*}\Bigr)\otimes \Q
\ee
where ${\rm Int}$ denotes the subspace of the integrable symbols. 

\bl
The shuffle product $\circ$ and the deconcatenation coproduct 
$\delta$
provide the graded space ${\bf I}_\bullet(X)$ with a structure of a commutative graded Hopf algebra. 
\el

{\bf Proof}. 
The space ${\bf I}_\bullet(X)$ is a subspace of the Hopf algebra 
$({\rm T}(\overline {{\cal O}_X^*})\otimes \Q, \circ, \delta)$. 
Clearly deconcatenation of an integrable symbol is an integrable symbol. 
It is easy to check that the shuffle product of integrable symbols is an integrable symbol. 
For example, given an integrable symbol $g\otimes h$, the shuffle product
$$
f\circ (g\otimes h) = f\otimes g\otimes h + g\otimes f\otimes h + g\otimes h\otimes f   
$$ is also integrable: 
projecting to $K_2$ modulo constants the first two factors of each summand we get zero since 
$\{f,g\} + \{g, f\} = 0$ and $\{g, h\}=0$ modulo constants by the assumption. 
The Lemma is proved. 

\paragraph{The Lie coalgebra of integrable symbols.} Let us consider the 
quotient of the graded commutative Hopf algebra ${\bf I}_\bullet(X)$ by the subspace 
${\bf I}_{>0}(X){\bf I}_{>0}(X)$ given by the products of the integrable systems of non-zero length:
$$
{\bf L}_\bullet(X):= \frac{{\bf I}_\bullet(X)}{{\bf I}_{>0}(X){\bf I}_{>0}(X)}.
$$ 
Then the coproduct on ${\bf I}_\bullet(X)$ determines a coproduct on the quotient, providing 
${\bf L}_\bullet(X)$ with a graded Lie coalgebra structure. We call it the Lie coalgebra of integrable symbols 
on $X$.

\subsection{Tate iterated integrals and variations of mixed Tate motives.} 

\paragraph{Variations of Hodge-Tate structures.}

Below we work with the category of $\Q$-Hodge-Tate variations. The key point is that it is a mixed Tate category, 
see Appendix in \cite{G2}. We also use the notion of 
the period of a  
variation of framed Hodge-Tate structures, see Section 4 of \cite{G6}. 
For convenience of the reader we recall now some of the basic properties. 
Below $X$ is a regular complex variety.

\vskip 2mm
A $\Q$-Hodge-Tate variation ${\cal V}$ is a variation of mixed Hodge $\Q$-structures. 
It  has a weight filtration denoted by $W_{\bullet}$. The  
associate graded ${\rm gr}^W_{-2m}{\cal V}$ are direct sums 
of the constant variations $\Q(m)_X$ of the rank one Hodge Tate structures of the Hodge type $(-m,-m)$ 
  on $X(\C)$, and ${\rm gr}^W_{-2m+1}{\cal V}=0$. 

An $n$-framing on a $\Q$-Hodge-Tate variation  ${\cal V}$ is a pair of non-zero morphisms 
$$
v: \Q(0)_X \lra {\rm gr}^W_{0}{\cal V}, \qquad f: {\rm gr}^W_{-2n}{\cal V} \lra \Q(n)_X.
$$
Let us consider the equivalence relation on the set of all $n$-framed $\Q$-Hodge-Tate variations
  on $X$ generated by the condition that a morphism 
of mixed Hodge structures  ${\cal V}_1 \to {\cal V}_2$ respecting the frames is an equivalence. 
Then the set of equivalence classes form a $\Q$-vector space, denoted by ${\cal H}_n(X)$. 
The addition is induced by the direct sum of variations. 
The tensor product induces a map ${\cal H}_n(X) \otimes {\cal H}_m(X) \lra {\cal H}_{n+m}(X)$, making 
$$
{\cal H}_\bullet(X) := \bigoplus_{n=0}^\infty{\cal H}_n(X) 
$$
into a commutative algebra over $\Q$, graded by the non-negative integers. Finally, there is a coassociative coproduct 
$$
\nu: {\cal H}_\bullet(X) \lra {\cal H}_\bullet(X)\otimes {\cal H}_\bullet(X)
$$
providing $({\cal H}_\bullet(X), \nu)$ with a structure of a commutative graded Hopf algebra. 
The category of graded comodules over this Hopf algebra is canonically equivalent to the category of 
$\Q$-Hodge-Tate variations on $X$. 

An $n$-framed $\Q$-Hodge-Tate variation   
  on  $X$ 
provides an element 
\be \la{I}
{\cal I} \in {\cal H}_n(X). 
\ee

\paragraph{Tate iterated integrals are periods of mixed Tate motives.} 

We say that a framed variation of mixed Hodge-Tate strucrtures is of {\it geometric origin} 
if it is equivalent to a one which can be realized in the cohomology of simplicial complex algebraic 
varieties. 

\bt \la{IIN} 
Let $X$ be a rational variety. Then,  given an integrable symbol 
$
{\rm I} \in {\bf I}_n(X), 
$ the Tate iterated integral $\int_a^b d\log ({\rm I})$ 
is a period of an $n$-framed  $\Q$-Hodge-Tate   
variation of geometric origin on $X\times X$:
\be \la{10.31.10.2}
{\cal I}({\rm I}) \in {\cal H}_n(X\times X).
\ee
This way we get an injective homomorphism of graded commutative algebras
$$
{\cal I}: {\bf I}_\bullet(X) \hra {\cal H}_\bullet(X\times X).
$$
\et

\vskip 3mm 
{\bf Warning}. The map ${\cal I}$ in general does not commute with the coproduct, even in the simplest case 
of the symbol $(t-a)\otimes (t-b)$ in $\Q(t)^*\otimes \Q(t)^*$. 

\vskip 3mm 

{\bf Proof}. The Tate iterated integral $\int_\gamma d\log {\rm I}$ is a period of the 
framed mixed Hodge structure provided by Beilinson's construction, see \cite{DG}. 
Namely, take the mixed Hodge structure 
${\cal P}(X;a,b)^*$ 
on the dual to the pronilpotent 
torsor of path between the base points $a,b$. The framing is given by 
the cohomology class 
$d\log ({\rm I})$ and the relative homology class provided by the homotopy 
class of a path $\gamma$ between $a$ and $b$. By construction, 
$\int_\gamma d\log ({\rm I})$ 
is the  period of a variation of framed mixed 
Hodge structures realized in the  cohomology of algebraic varieties. We do not claim 
however that  the mixed Hodge structure ${\cal P}(X;a,b)^*$ is Hodge-Tate. 
We claim only that it is equivalent to a Hodge-Tate one. This implies that there is a canonical minimal  
(see, say, the Appendix to \cite{G2}) 
Hodge-Tate representative in the equivalence class, providing a Hodge-Tate variation on
 $X\times X$, uniquely defined by the symbol. The latter 
is the same thing as an element (\ref{10.31.10.2}). 

So let us show that for any $(a,b)\in X\times X$ 
the obtained $n$-framed mixed Hodge structure is equivalent to a Hodge-Tate one.
Since $X$ is rational, 
there exists a punctured rational curve $C$ on $X$ connecting the points $a,b$. 
The mixed Hodge structure on ${\cal P}(C;a,b)^*$ is evidently Hodge-Tate.  
The canonical morphism ${\cal P}(X;a,b)^* \lra {\cal P}(C;a,b)^*$ induces an equivalence of framed 
Hodge-Tate structures. The injectivity of the map ${\cal I}$ is obvious. The claim that it is a homomorphism 
of algebras 
follows from the product formula for the framed Hodge-Tate variations assigned to the iterated integrals 
on the line \cite{G2}. 
The theorem is proved.

\bcon \la{MCONJ} For any variety $X$, the Tate iterated integral $\int_x^y d\log ({\rm I})$ 
is the  period of a framed $\Q$-Hodge-Tate variation  
on  $X\times X$. 
\econ

{\bf Remark}. The argument above shows that Conjecture \ref{MCONJ} reduces to the case when $X$ is a curve. 
For the length   one iterated integral, i.e. $n=1$ in (\ref{IEL}),  it is obvious, and for $n=2$ it is 
easy to prove. So $n=3$ is the first non-trivial case.

\vskip 3mm
Thanks to Theorem \ref{IIN}, an integrable symbol ${\rm I}$ on $X$ gives rise to an 
$n$-framed geometric Hodge-Tate variation   ${\cal I}({\rm I})$ on  
$X \times X$. We denote by ${\cal I}_{x,y}({\rm I})$ its fiber at the point $(x,y)$. 
Let us show how to  recover  the symbol ${\rm I}$ 
from the $n$-framed variation ${\cal I}({\rm I})$.

We start from a general construction assigning a  {\it symbol} on $X$ 
to any $n$-framed Hodge-Tate variation on $X$. Then we show that  
the {\it reduced  symbol} of the variation ${\cal I}_{x,y}({\rm I})$ on 
$X \times X$ does not 
depend on the first factor, and recovers 
the original symbol ${\rm I}$.

\paragraph{Symbol of a framed Hodge-Tate variation.} 
Recall (Appendix in \cite{G2}) that 
$$
{\cal H}_0(X) = \Q, \qquad {\cal H}_1(X) = {\cal O}^*_{X, {\rm an}}\otimes \Q.
$$
Here ${\cal O}^*_{X, {\rm an}}$ 
denotes the multiplicative group of analytic functions on $X$. So the iterated coproduct provides us a map
$$
\nu_{[n]}: {\cal H}_n(X) \lra \bigotimes^n {\cal O}^*_{X, {\rm an}}\otimes \Q. 
$$
\bd
The symbol $S_n({\cal V})$ of an $n$-framed  Hodge-Tate variation 
${\cal V} $ on $X$ is given by 
$$
S_n({\cal V}):= \nu_{[n]}{\cal V} \in \bigotimes^n {\cal O}^*_{X, {\rm an}}\otimes \Q.
$$
\ed

Here is a  way to calculate the symbol. 
An $(m-1, m)$-framing 
on ${\cal V}$ is a pair of non-zero maps
\be \la{framing}
e: \Q(m-1)_X \lra {\rm gr}^W_{-2m+2}{\cal V}, \qquad 
f: {\rm gr}^W_{-2m}{\cal V}\lra \Q(m)_X.
\ee
 It gives rise to an extension  class 
$$
e(s,f) \in {\rm Ext}_{\Q-MHS}^1(\Q(0)_X, \Q(1)_X) \stackrel{\sim}{=} {\cal O}^*_{X, {\rm an}}\otimes \Q. 
$$

Given an $n$-framed Hodge-Tate variation  ${\cal V}$, choose a basis $\{e^{(k)}_{\bullet}\}$ 
in ${\rm gr}^W_{-2k}{\cal V}$ for each $-n\leq k\leq 0$, so that for $k=0$ and $k=-n$ it 
coincides with the given framing. Let $\{f^{(k)}_{\bullet}\}$ be the dual basis. Then 
\be 
S_n({\cal V}) = \sum_i\bigotimes_{k=-n}^0 e(e^{(k)}_{i}, f^{(k)}_{i})
\ee
where the sum is over all basis elements. 
The proof follows easily by induction by applying the coproduct to ${\cal V}$.

\paragraph{The reduced Hopf algebra $\overline {\cal H}_{\bullet}(X)$.} 
Let us set 
$$
\overline {{\cal O}_{X, {\rm an}}^*}:=
{\cal O}_{X, {\rm an}}^*/\C^*.
$$ 

\bd Let $X$ be a regular complex variety. The algebra 
$\overline {\cal H}_{\bullet}(X)$ is the quotient of  the algebra ${\cal H}_{\bullet}(X)$ 
by the ideal generated by constant variations -- the latter is canonically isomorphic to ${\cal H}_{\bullet}({\rm Spec}(\C))$, by restriction to any point of $X$. 
\ed

Clearly the algebra 
$\overline {\cal H}_{\bullet}(X)$  is a Hopf algebra.

Given a point $a\in X$, and varying a point $z\in X$, the $n$-framed motivic iterated integrals provide an element 
\be \la{10.31.10.3a}
{\cal I}_{a,z} (f_1 \otimes f_2 \otimes  \ldots \otimes f_n) \in {\cal H}_n(X).
\ee

\bl \la{10.31.10.100}
The image 
$$
\overline {\cal I}_{a,z} (f_1 \otimes f_2 \otimes  \ldots \otimes f_n) \in \overline {\cal H}_{n}(X)
$$
of the $n$-framed variation (\ref{10.31.10.3a})  does not depend on the 
choice of the point $a$. 
\el

{\bf Proof}. There is a formula for motivic iterated integrals, understood as elements of 
${\cal H}_{\bullet}$, where $a,b,z$ are any points in $X$:
\be \la{10.31.10.3}
{\cal I}_{a,z} (f_1 \otimes f_2 \otimes  \ldots \otimes f_n) = 
\sum_{k=0}^n
{\cal I}_{a,b}  (f_1 \otimes f_2 \otimes  \ldots \otimes f_k)
 \cdot {\cal I}_{b,z} (f_{k+1} \otimes \ldots \otimes f_n). 
\ee
The Lemma follows from this formula. Indeed, the ${\cal I}_{a,b}$ here is a constant, so all summands with $k>0$ die in 
$\overline {\cal H}_{n}(X)$. 

\bd 
The reduced symbol $\overline S_n({\cal V})$ of an $n$-framed Hodge-Tate variation ${\cal V}$ 
is the projection of the symbol $S_n({\cal V})$ to 
$\bigotimes^n\overline {{\cal O}_{X, {\rm an}}^*}$.
\ed

The reduced symbol $\overline S_n({\cal V})$ is nothing else but the iterated coproduct 
$\nu_{[n]}$ applied to the image $\overline {\cal V} \in \overline {\cal H}_n$ 
of ${\cal V}\in {\cal H}_n$.

\vskip 3mm
Thanks to Lemma \ref{10.31.10.100} the projection of the motivic iterated integral 
${\cal I}_{a,z}({\rm I})$ to the reduced 
Hopf algebra $\overline {\cal H}_n$ is independent of $a$, and provides a homomorphism of abelian groups 
$$
\overline {\cal I}_n: {\bf I}_n(X) \lra \overline {\cal H}_n(X).
$$
 On the other hand, the reduced symbol is a homomorphism of abelian groups 
$$
\overline S_n: \overline {\cal H}_n(X) \lra \bigotimes^n\overline {{\cal O}_{\rm an}^*}(X)\otimes \Q.
$$

\bt \la{VITJ}
Assume that $X$ is rational. Then $\overline S_n \circ \overline {\cal I}_n$ is the identity map. 
The map 
$$
\overline {\cal I}= \oplus\overline {\cal I}_n: {\bf I}_\bullet(X) 
\stackrel{}{\lra} \overline {\cal H}_\bullet(X) 
$$ is an injective homomorphism of Hopf algebras. 
\et

{\bf Proof}. Let us assume that $X$ is a punctured projective line. The main result of \cite{G2} 
describes the coproduct and therefore the symbol of the motivic iterated integrals on the line. 
The claim that the composition $\overline S_n \circ \overline {\cal I}_n$ is the identity map, as well as the claim that 
 $\overline {\cal I}$ is a homomorphism of Hopf algebra follows 
immediately from this. The general case is reduced to the case of the  
punctured projective line, since a symbol is determined by its restriction to the generic projective line in $X$.  
Set  $\overline {S}:= \sum_n\overline {S}_n$. Since the composition
$$
{\bf I}_\bullet(X) 
\stackrel{\overline {\cal I}}{\lra} \overline {\cal H}_\bullet(X) \stackrel{\overline {S}}{\lra} {\bf I}_\bullet(X)
$$
is the identity map, the map $\overline {\cal I}$ is injective. 
Theorem \ref{VITJ} is proved.

\paragraph{Weakly geometric Hodge-Tate variations.} There is a natural map
\be \la{IMM}
{\cal O}_X^*\otimes \Q \hra {\rm Ext}_{\Q-MHS}^1(\Q(0)_X, \Q(1)_X)
\ee
where the Ext group is in the category of variations 
of mixed $\Q$-Hodge structures on $X(\C)$.  

\bd \la{GEOVAR} A framed $\Q$-Hodge-Tate variation  on $X(\C)$ is weakly geometric if 
the ${\rm Ext}^1$ 
defined by any $(m-1, m)$-framing is in the image of map (\ref{IMM}). 
\ed

Denote by ${\cal H}^{\rm wg}_\bullet(X)$ 
the Tannakian Hopf algebra of the category of weakly geometric 
$\Q$-Hodge-Tate variations on  $X$. 
One has (Appendix in \cite{G2})
$$
{\cal H}^{\rm wg}_0(X) = \Q, \qquad {\cal H}^{\rm wg}_1(X) = {\cal O}^*_X\otimes \Q.
$$
So the iterated coproduct is a map
$$
\nu_{[n]}: {\cal H}^{\rm wg}_n(X) \lra \bigotimes^n {\cal O}^*_X\otimes \Q. 
$$

{\bf Remarks}. 1. The map (\ref{IMM}) should provide an isomorphism
\be \la{IMM1}
{\cal O}_X^*\otimes \Q \stackrel{\sim}{\lra} {\rm Ext}_{\Q-Mot}^1(\Q(0)_X, \Q(1)_X)
\ee
 where on the right hand side we have the {\rm Ext}-group in the (say, Voevodsky) category of 
mixed motivic sheaves. However, although we have such a map, and it is 
injective,  we do not know its surjectivity.

2. There are constant variations of Hodge-Tate structures over a regular complex variety $X$ which are not motivic. 
For example, ${\rm Ext}_{\Q-MHS}^1(\Q(0), \Q(2))=\C/(2\pi i)^2\Q$, while the Hodge realization 
of ${\rm Ext}_{\Q-Mot}^1(\Q(0), \Q(2))$ 
is smaller, countable, due to the rigidity of the regulator map.

\paragraph{A conjectural description of the Hopf algebras $\overline {\cal H}^{\rm g}_\bullet$ and 
$\overline {\cal H}^{\rm wg}_\bullet$.} 

Denote by ${\cal H}^{\rm g}_n(X)$ the Tannakian Hopf algebra 
of variations of Hodge-Tate structures of geometric origin. Clearly there is an inclusion 
$i: {\cal H}^{\rm g}_n(X) \hra {\cal H}^{\rm wg}_n(X)$.

\bcon
The inclusion $i$ gives rise to an isomorphism  $\overline i: \overline 
{\cal H}^{\rm g}_n(X) \stackrel{\sim}{\lra} \overline {\cal H}^{\rm wg}_n(X)$.  
\econ

\bcon \la{13:42}
The sum of the maps $\overline S_n$ provides an isomorphism 
\be \la{13:44}
\overline S: \overline {\cal H}^{\rm wg}_\bullet(X) \stackrel{\sim}{\lra} {\bf I}_\bullet(X).
\ee 
\econ

Notice that we do not know that the map $d\log: K_2(\Q(X)) \lra \Omega^2_{\rm log}(X)$ is injective 
even for $X={\Bbb A}^2$. So we can not prove that the image of the map $\overline S$ consists of 
integrable elements.  

\vskip 3mm
{\bf Remark}. By Lemma \ref{11.1.10.5} the map 
(\ref{13:44}) is a morphism of Hopf algebras. So Conjecture \ref{13:42} tells that the map 
(\ref{13:44}) should be an isomorphism of Hopf algebras.

\paragraph{Differential equation for the period.} Let 
$p({\cal V})$ be the multivalued analytic function 
on $X(\C)$ given by the period of a framed variation 
 ${\cal V}$. The period functions assigned to equivalent variations 
are the same. 
Therefore there is a map $p \otimes d\log$ from 
${\cal H}_{n-1}(X) \otimes {\cal O}_X^*$ to multivalued analytic 
$1$-forms at the generic point of $X(\C)$. 
Let 
$$
\nu_{n-1, 1}: {\cal H}_{n}(X) \lra  {\cal H}_{n-1}(X) \otimes {\cal O}_X^*. 
$$
be the $(n-1, 1)$-component of the coproduct. The following is Lemma 4.6a) in \cite{G6}.

\bl \la{3} The differential of the period $p({\cal I})$ of a framed Hodge-Tate variation 
${\cal I}$ is given by  
$$
dp({\cal I}) = p\otimes d\log\Bigl( \nu_{n-1, 1}({\cal I})\Bigr). 
$$
\el

\section{The symbol of the Grasmannian polylogarithm}

\subsection{The symbol of the Grassmannian $n$-logarithm}
Recall that a point of the Grassmannian $G^n_n$ can be described as a configuration $(l_1, ... , l_{2n})$ of $2n$ 
vectors in an $n$-dimensional complex vector space.  
\bd A symbol
$
{\rm I}_n(l_1, ... , l_{2n}) \in \bigotimes^n {\cal O}(G^n_n)^*
$ 
is given by the formula
\be \la{MF}
{\rm I}_n(l_1, ...  l_{2n}):= 
{\rm Alt}_{2n}\Bigl(\Delta(l_1, ... ,l_{n-1}, l_n) \otimes 
\Delta(l_2, ... , l_{n+1}) \otimes ... \otimes 
\Delta(l_{n}, ... , l_{2n-1})\Bigr).
\ee
\ed

\paragraph{Comparison Theorem.} 
It relates $\Lambda_{n-1,1}$ and ${\rm I}_n$. Observe that 
it is sufficient to know $\nu_{n-1,1}$ in order to compute $\nu_{[n]}$. Indeed, $\nu_{[n]}$ is the composition 
$\ldots \circ (\nu_{n-3,1}\otimes {\rm Id}\otimes {\rm Id})\circ (\nu_{n-2,1}\otimes {\rm Id})\circ \nu_{n-1,1}$.
\bt \label{keyl}
One has 
$$
(\nu_{[n-1]} \otimes {\rm Id} ) \circ \Lambda_{n-1,1} (l_1, ... ,l_n , m_1, ... , m_n) ~= ~
(-1)^{n}2(n!)^2 {\rm I}_n(l_1, ... , l_n; m_1, ... ,  m_n). 
$$
\et

{\bf Proof}. Using (8) and (10) to calculate $\nu_{n-2,1}$, continuing the same line,  
to calculate $\nu_{n-3,1}$ of the first factor, end so on, we get the following expression for 
the term in $A_2 \otimes F^* \otimes ... \otimes F^*$: 
$$
(-1)^{n-3}n^2 \cdot ... \cdot 4^2 \cdot {\rm Alt}_{2n}
\Bigl(\langle l_1,... , l_{n-3}|l_{n-2}, l_{n-1}, l_n; m_{n-2}, m_{n-1}, m_{n}\rangle _{A_{2}}  \otimes 
$$
$$
\Delta (l_1, ..., l_{n-3}, m_{n-2}, m_{n-1} ,m_{n})\otimes  ... \otimes \Delta (m_{1}, ... ,m_{n})\Bigr).
$$
Taking into account formula (\ref{BNU}) for $  \nu_{1,1}$, with the footnote 2), we get 
\begin{equation} \label{EQQ1}
(-1)^{n-2}n^2 \cdot ... \cdot 4^2\cdot 3^2 {\rm Alt}_{2n}
\Bigl(
\Delta (l_1, ..., l_{n-3}, m_{n-2}, l_{n-1}, l_{n}) \otimes 
\langle l_1,... , l_{n-3}, m_{n-2}|l_{n-1}, l_n; m_{n-1}, m_{n}\rangle _{A_{1}}   
\end{equation}
$$
\otimes
\Delta (l_1, ..., l_{n-3}, m_{n-2}, m_{n-1},m_{n})\otimes ... \otimes \Delta (m_1, m_{2}, ... ,m_{n}) +
$$
\begin{equation} \label{EQQ2}
\langle l_1,... , l_{n-2}|l_{n-1}, l_n; m_{n-1}, m_{n}\rangle _{A_{1}}  \otimes 
\Delta (l_1, ..., l_{n-2}, m_{n-1}, m_{n})\otimes ... \otimes \Delta (m_1, m_{2}, ... ,m_{n})\Bigr). 
\end{equation}
Using the formula
\begin{equation} \label{EQQ3}
\langle l_1,..., l_{n-2}|l_{n-1}, l_n; m_{n-1}, m_{n}\rangle_{A_1} = 
\frac{\Delta (l_1, ..., l_{n-2}, l_{n-1}, m_{n-1})\Delta (l_1, ..., l_{n-2}, l_{n}, m_n)}{\Delta (l_1, ..., l_{n-2}, l_{n-1}, m_n)
\Delta (l_1, ..., l_{n-2}, l_{n}, m_{n-1})}
\end{equation}
we write the term (\ref{EQQ1}) as follows
$$
-(-1)^n (n!)^2{\rm Alt}_{2n}\Bigl(
\Delta (l_1, ..., l_{n-3}, m_{n-2}, l_{n-1}, l_{n})\otimes \Delta (l_1, ..., l_{n-3}, m_{n-2}, l_{n-1}, m_{n}) 
\otimes 
$$
$$
\Delta (l_1, ..., l_{n-3}, m_{n-2}, m_{n-1}, m_{n}) \otimes ... \otimes 
\Delta (m_1, m_{2}, ... ,m_{n})\Bigr) = 
$$
\begin{equation} \label{EQQ4}
-(-1)^n 2(n!)^2{\rm Alt}_{2n} \Bigl(
\Delta (l_1, ..., l_{n-1}, m_{n})\otimes ... 
\otimes 
\Delta (m_1, m_{2}, ... ,m_{n})\Bigr). 
\end{equation}
In the last step we use the fact that each of the permutations $(l_{n-2}, l_{n-1}, l_{n}) \lra (l_{n}, l_{n-2}, l_{n-1})$ and 
$(m_{n-2}, m_{n-1}, m_{n}) \lra (m_{n}, m_{n-2}, m_{n-1})$ are even. 
Theorem \ref{keyl} is proved.

\bt \la{10.31.10.1}
a) The symbol ${\rm I}_n$ is integrable. 

b) It lives on  
$PG_n$, and satisfies two 
$(2n+1)$-term relations: 

\noindent
1) For a generic 
configuration of $2n+1$ vectors $(l_1,...,l_{2n+1})$ in 
$V_{n}$ one has 
\be \la{1}
\sum_{i=1}^{2n+1}(-1)^i {\rm I}_n(l_1,...,\widehat l_i, ... , l_{2n+1}) = 0.
\ee
2) For a  generic configuration 
of vectors $(m_1,...,m_{2n+1})$ in $V_{n+1}$ one has 
\be \la{2}
\sum_{j=1}^{2n+1}(-1)^j {\rm I}_n(m_j| m_1,...,\widehat m_j, ... , m_{2n+1}) = 0.
\ee
\et

{\bf Proof}.  a) Follows easily from Lemma \ref{MMLEa} by using 
Comparison Theorem \ref{keyl}. 

\vskip 2mm
b) Changing
 the vector $l_1$ to $al_1$ we get  
$$
{\rm I}_n(al_1, ...  l_{2n}) - {\rm I}_n(l_1, ...  l_{2n}) = 
{\rm Alt}_{2n}\Bigl(a \otimes 
\Delta(l_2, ... , l_{n+1}) \otimes ... \otimes 
\Delta(l_{n}, ... , l_{2n-1})\Bigr) = 0.
$$
Indeed, we skewsymmetrize an expression which does not 
contain the pair of vectors $(l_1, l_{2n})$.  

The two relations follow immediately from 
Comparison Theorem \ref{keyl} and Lemma \ref{keyp1}. 
Theorem \ref{10.31.10.1} is proved.

\paragraph{Conclusion.} The iterated integral assigned to the symbol ${\rm I}_n$ is a 
multivalued analytic function at the generic point of  $G_n(\C) \times G_n(\C)$. By Theorem \ref{IIN} 
it is the period of a motivic variation of framed Hodge-Tate structures at the generic point of 
$G_n(\C) \times G_n(\C)$. Modulo the ideal of constant variations, 
it is a variation at the generic point of $G_n(\C)$. 

\subsection{The symbol of the bi-Grassmannian $n$-logarithm cocycle}

We conjecture that there exists a nice explicit expression for the symbol 
of the bi-Grassmannian $n$-logarithm cocycle. Let us formulate this precisely. 

Recall the Lie coalgebra ${\bf L}_\bullet(X)$ of integrable symbols on a variety $X$. 
There is the standard cochain complex of the Lie coalgebra ${\bf L}_\bullet(X)$:
$$
{\bf L}_\bullet(X) \stackrel{\delta}{\lra} \Lambda^2 {\bf L}_\bullet(X)\stackrel{\delta}{\lra}
 \Lambda^3 {\bf L}_\bullet(X)\stackrel{\delta}{\lra} \ldots
$$
Here the first map is the coproduct, and the other maps are induced by the coproduct via the Leibniz rule. 

Recall the 
Grassmannian  $G_p^q$ of $q$-dimensional subspaces in a coordinate vector space of 
dimension $p+q$, transversal to the coordinate hyperplanes. There are maps between the Grassmannians
$$
A_i: G_p^q \lra G_p^{q-1}, \qquad B_j: G_p^q \lra G_{p-1}^{q}, \quad 1 \leq i,j \leq p+q.
$$
The map $A_i$ is given by the intersection with the $i$-th coordinate hyperplane, and $B_j$ is 
ionduced by the projection 
along the $j$-th coordinate axis. 

Recall that a point of the Grassmannian $G_p^q$ can be encoded by a configuration of $p+q$ vectors 
$(l_1, ..., l_{p+q})$ in a vector space of dimension $p$. Then one has 
$$
A_i(l_1, ..., l_{p+q}) = (l_1, ..., \widehat l_i, ... l_{p+q}), \qquad 
B_j(l_1, ..., l_{p+q}) = (l_j ~|~ l_1, ..., \widehat l_j, ... l_{p+q})
$$
\bcon
Given a positive integer $n$, there exist elements of total weght $n$
$$
{\rm I}(n)^q_{p} \in \Lambda^{2n+1 - p -q }{\bf L}_\bullet({\rm G_p^q}), \qquad {\rm I}(n)^q_{p} =0 ~~\mbox{for $p<n$},
$$
satisfying the following conditions:

\begin{itemize}
\item  The bi-Grassmannian cocycle condition (here $A_i^*$ and $B_j^*$ are the pull backs):
$$
\delta {\rm I}(n)^q_{p} = \sum_{i=1}^{p+q}(-1)^i A_i^*{\rm I}(n)^{q-1}_{p} + \sum_{j=1}^{p+q}(-1)^j B_j^*{\rm I}(n)^{q}_{p-1}. 
$$

\item The symbol ${\rm I}(n)^n_{n} \in {\bf L}_n({\rm G_n^n})$ is the 
symbol of the Grassmannian $n$-logarithm function: 
\be \la{MF1}
{\rm I}(n)^n_n(l_1, ...  l_{2n}):= 
{\rm Alt}_{2n}\Bigl(\Delta(l_1, ... , l_n) \otimes 
\Delta(l_2, ... , l_{n+1}) \otimes ... \otimes 
\Delta(l_{n}, ... , l_{2n-1})\Bigr).
\ee

\item The symbol ${\rm I}(n)^1_{n} \in \Lambda^n{\bf L}_n({\rm G_n^1})$ is given by the formula
$$
{\rm I}(n)^1_{n}(l_1, ..., l_{n+1}) = {\rm Alt}_{n+1} \Bigl(\Lambda_{i=1}^{n}\Delta(l_1, ..., \widehat l_i, ..., l_{n+1})\Bigr).
$$
\end{itemize} 
\econ

A construction of these cocycles for $n\leq 4$  folows from the main results of \cite{G} and \cite{G3}.


\begin{thebibliography}{BGSV}

\bibitem[A]{A} Aomoto K.: {\it Addition theorem of Abel type for hyper-logarithms}.  
Nagoya Math. J.  88  (1982), 55--71. 

\bibitem[BMS]{BMS} Beilinson A.A., MacPherson R., Schechtman, 
V.V: {\it Notes on motivic cohomology}, Duke Math.\ J.\ 54 (1987), 
679--710. MR0899412 (88f:14021).

 \bibitem[BGSV]{BGSV} Beilinson A.A., Goncharov A.B., Schechtman, 
V.V, Varchenko A.N.: {\it Aomoto dilogarithms, mixed Hodge structures and motivic cohomology of pairs of triangles on the plane}.  
The Grothendieck Festschrift, Vol. I,  135--172, Progr. Math., 86, Birkhäuser Boston, Boston, MA, 1990. 

 \bibitem[Ch]{Ch} Chen K-T.: {\it Algebras of iterated path integrals and 
fundamental groups}. Trans. Amer. Math. Soc. 156(1971), 359--379. 


\bibitem[DG]{DG} Deligne P., Goncharov A.B.: {\it 
Groupes fondamentaux motiviques de Tate mixte.}
  Ann. Sci. \'Ecole Norm. Sup. (4)  38  (2005),  no. 1, 1--56. 

\bibitem[GGL]{GGL} Gabrielov A., Gelfand I.M., Losik M.: 
{\it Combinatorial computation of characteristic classes} I, II,
Funct. Anal. i ego pril. 9 (1975) 12-28, ibid 9 (1975) 5-26. MR0410758 (53 14504a). 



\bibitem[GM]{GM} Gelfand I.M., MacPherson R: {\it Geometry in 
Grassmannians and a generalisation of the dilogarithm}, Adv. 
in Math., 44 (1982) 279--312. MR0658730 (84b:57014).



\bibitem[G]{G} Goncharov A.B.: {\it Geometry of configurations, polylogarithms, and motivic cohomology}.
Adv. Math. 114 (1995), no. 2, 197â318. 

\bibitem[G1]{G1} Goncharov A.B.: {\it Chow polylogarithms and
regulators}. Math. Res. Letters, 2, (1995), 99-114. MR1312980 (96b:19007). 

\bibitem[G2]{G2} Goncharov A.B.: {\it Galois symmetries of fundamental 
groupoids and noncommutative geometry}.  Duke Math. J.  128  (2005),  no. 2, 209--284. 

\bibitem[G3]{G3} Goncharov, A. B. {\it Geometry of the trilogarithm and the motivic Lie algebra of a field}.  Regulators in analysis, geometry and number theory,  127--165, 
Progr. Math., 171, Birkhäuser Boston, Boston, MA, 2000.

\bibitem[G4]{G4} Goncharov, A. B. 
{\it Polylogarithms, regulators and Arakelov motivic complexes.} JAMS  J. Amer. Math. Soc.  18  (2005),  no. 1, 1--60. arXiv:math/0207036.

\bibitem[G5]{G5} Goncharov, A. B. {\it Explicit construction of characteristic classes}. 
{Advances in Soviet Mathematics,} 16, v.\ 1, 
    Special volume dedicated to 
    I.M.Gelfand's 80th birthday, 169 - 210, 1993.

\bibitem[G6]{G6} Goncharov, A. B. {\it Volumes of hyperbolic manifolds and mixed Tate motives} 
JAMS vol. 12, N2, (1999), 569-618. arXiv:alg-geom/9601021.


\bibitem[HM]{HM} Hain R, MacPherson R: {\it Higher Logarithms}, 
Ill. J.\ of Math,, vol. 34, (1990) N2,  392--475. MR1046570 (92c:14016). 

\bibitem[H]{H} Hain R: {\it The existence of higher logarithms}. 
Compositio Math. 100 (1996) N3, 247- 276.  
MR1387666 (97g:19008).

\bibitem[HY]{HY}  Hain, R. Yang, J.: {\it Real Grassmann polylogarithms 
and Chern classes}.  Math. Ann.  304  (1996),  
no. 1, 157--201.

\bibitem[HM1]{HM1} Hanamura M, MacPherson R: {\it Geometric construction
of polylogarithms}, Duke Math. J. 70 ( 1993)481-516. MR1224097 (94h:19005). 




\bibitem[HM2]{HM2} Hanamura M,  MacPherson R.: 
{\it Geometric construction of polylogarithms, II}, 
Progress in Math. vol. 132 (1996) 215-282. MR1389020 (97g:19007). 



\bibitem[M]{M} MacPherson R.: {\it Gabrielov, Gelfand, Losic
combinatorial formula for the first Pontryagin class} Seminar Bourbaki 1976. MR0521763 (81a:57022).

\bibitem[Su]{Su} Suslin A.A: {\it Homology of $GL_{n}$,
characteristic classes and Milnor's $K$-theory.} Trudy Mat. Inst. Steklov, 
165 (1984), 188-204;

\bibitem[Su2]{Su2} Suslin A.A: {\it K-theory of a field and the Bloch group}, 
Proc. Steklov Inst. Math. 4 (1991) 217 239. 

\bibitem[Yu]{Yu} Yusin, Boris V.: {\it  Sur les formes $S\sp{p,q}$ apparaissant dans le calcul combinatoire de la deuxi\'eme classe de Pontriaguine par la m\'ethode de 
Gabrielov, Gelfand et Losik}. (French)  C. R. Acad. Sci. Paris Sér. I Math.  292  (1981), no. 13, 641--644. 57R20



\end{thebibliography}
\end{document}